\documentclass{amsart}

\usepackage{enumerate}
\usepackage{amsfonts}
\usepackage{amssymb, amscd}
\usepackage{graphicx}

\usepackage{float}

\newcommand{\pf}{\n{\em Proof.}   }

\newcommand{\n}{\noindent}

\newcommand{\sn}{\smallskip\noindent}

\newtheorem{teo}{Theorem}[section]
\newtheorem{cor}[teo]{Corollary}
\newtheorem{prop}[teo]{Proposition}
\newtheorem{lema}[teo]{Lemma}

\theoremstyle{definition}

\title{Complex ellipsoids and complex symmetry}

\author{Jorge Arocha}
\author{Javier Bracho}
\author{Luis Montejano}

\begin{document}
\maketitle
\begin{abstract}
Several characterizations of complex ellipsoids among convex bodies in $\mathbb{C}^n$, in terms of their sections and projections are proved. Characterizing complex symmetry in similar terms is an important tool.
\end{abstract}

\section{Introduction}

An ellipsoid is the image of a ball under an affine transformation. If this affine transformation is over the complex numbers, we refer to it as a \emph{complex ellipsoid}. Characterizations of  real ellipsoids have received much attention over the years (see Soltan's survey, \cite{So}); however, characterizations of complex ellipsoids have been scarcely studied. In the paper of Bracho and Montejano \cite{BM} on the proof of the Complex Banach  Conjecture, a characterization of complex ellipsoids plays an important role. Other antecedents are Gromov's paper \cite {G} and our recent paper on extremal inscribed and circumscribed complex ellipsoids \cite{ABM}. 

Our main contribution is an unexpected characterization of complex ellipsoids
which has no analogue over the real numbers. Namely, we prove that if every complex line intersects a convex body  in complex space $\mathbb{C}^n$ in a disk, then the convex body is a complex ellipsoid (Theorem~\ref{teoBom}).

The study of complex ellipsoids is naturally related to the study of complex symmetry. So, characterizing and understanding complex symmetry is vital to characterizing complex ellipsoids. In Section~2 we give conditions to know when a convex body is complex symmetric. In particular, we prove one result in terms of sections (Theorem~\ref{lemBB}), and one in terms of projections (Theorem~\ref{teoSim}). 

In \cite{BM} it was proved that if every hyperplane section through the center of a complex symmetric convex body in complex space $\mathbb{C}^n$ is a complex ellipsoid, then the body is a complex ellipsoid. The proof relies heavily on the notion of complex symmetry.  In Secction~4, we use our main characterization to drop this condition, and prove the result in its full generality (Theorem~\ref{teoksec}). Also, we prove the analogue characterization in terms of orthogonal projections 
(Theorem~\ref{teokproj}), and another one in terms of affine equivalence of complex line sections (Theorem~\ref{entenado}). 
   
It is interesting to note the use of various types of topological techniques in the proofs of most of the results of this paper.

From now on and unless otherwise stated, the notions of  \emph{line},  \emph{plane},   \emph{hyperplane},  \emph{subspace}, 
 \emph{ellipsoid},  \emph{linear map} and  \emph{affine map}, are understood as complex.  Furthermore, a  \emph{convex body} is a compact, convex set with nonempty interior contained in complex space $\mathbb{C}^n$; and a  \emph{$k$-section} of it, is its intersection with an affine subspace of dimension $k$.

\section{Symmetry by means of sections and of projections}

 Let $\mathbb{S}^1$ be the space of all unit complex numbers. 
Let $A\subset \mathbb{C}^n$ be a set.  We say that $A$ is \emph{symmetric} if there is a translated copy $A'$ of $A$ such that 
$\xi A'=A'$, for every $\xi\in \mathbb{S}^1$. In this case, if  $A'=A-x_0$, we say that $x_0$ is the \emph{center of symmetry} of $A$.
If $-A$ is a translated copy of $A$, then we just say that $A$ is \emph{real-symmetric}.    It will be useful to consider the empty set as a symmetric set.  Note that \emph{a convex set $A\subset \mathbb{C}^n$ is symmetric with center at $x_0$  if and only if for every line $L$ trough $x_0$, the section $L\cap A$ is a disk centered at $x_0$.}  Of course, any hyperplane or  the unit ball of a finite dimensional Banach space over the complex numbers $\mathbb{C}$ is symmetric. 

We shall first prove that if every hyperplane section through a point of a convex body in complex space is symmetric, then the body is also symmetric (Theorem~\ref{lemBB}). For which we need the following useful characterization of symmetry. 

\begin{teo}\label{teosc}
A convex body $K$ is symmetric if and only if for every $\xi\in \mathbb{S}^1$, $\xi K$ is a translated copy of $K.$
\end{teo}

\pf Suppose $K$ is  symmetric, then there exists $x_0\in\mathbb{C}^n$ such that for every $\xi\in \mathbb{S}^1$, $\xi(K-x_0)=K-x_0$. This implies that 
$\xi K-\xi x_0=K-x_0$ and hence that 
$$\xi K= K + (\xi -1)x_0\,.$$ 
Thus, $\xi K$ is a translated copy of $K$. 
  
Suppose now, that for every $\xi\in \mathbb{S}^1$, $\xi K$ is a translated copy of $K$. Hence 
$$\xi K= K+ \kappa_\xi\,,$$ 
where $ \kappa_\xi\in \mathbb{C}^n$  depends on $\xi$. 

For $\xi = -1$, we get that $- K= K+ \kappa_{-1}$. This implies that $K$ is real-symmetric with center $-\frac{\kappa_{-1}}{2}$. Define $x_0=-\frac{\kappa_{-1}}{2}$; it is the center of real-symmetry of $K$.

We then have that for any given $\xi\in \mathbb{S}^1$, $\xi K$  is real-symmetric with center $\xi x_0$.
On the other hand, since $\xi K= K+ \kappa_\xi$, the center of real-symmetry of $\xi K$ is $x_0 + \kappa_\xi$. Since there is a unique center of real-symmetry of $\xi K$, we have that
$$\xi x_0=x_0 + \kappa_\xi\,.$$  
Therefore, $\kappa_\xi=(\xi -1)x_0.$

To conclude the proof, note that $\xi K=K+(\xi-1)x_0,$ implies that $\xi( K-x_0)=K-x_0$ for every $\xi\in \mathbb{S}^1$, and therefore that 
$K$ is symmetric, by definition. \qed

The following theorem is the complex version of Roger's Theorem \cite{R}.

\begin{teo}[Complex Roger's Theorem] \label{lemrg}
Let $K_1, K_2\subset \mathbb{C}^n$, $n\geq 3$, be strictly convex bodies with ${p_1\in int K_1}$ and $p_2\in int K_2$. Suppose that for every  hyperplane $H$ through the origin, the section $(H+p_1)\cap K_1$ is positively homothetic to $(H+p_2)\cap K_2$. Then, $K_1$ is positively homothetic of $K_2.$ 
\end{teo}

\pf 
Note that if $K_2'$ is positively homothetic to $K_2$, then $K_2'$ also satisfies the conditions of the lemma. Then, using a positive homothecy, we may assume that $p_1=p_2$ and that  $K_2\subset$ int$K_1$.  Now, dilate $K_2$ until it touches the boundary of $K_1$. So, without loss of generality we may assume that:
\begin{enumerate}
\item $K_2\subset K_1,$
\item $p_1=p_2$, and 
\item $K_1$ and $K_2$ share a common support real-hyperplane $\Delta$ at $q$. 
\end{enumerate}

Let $L$ be the line trough $q$ and $p_1=p_2$. By hypothesis, for every hyperplane $H$ through $L$ there is a positive homothecy  that sends $H\cap K_1$ onto $H\cap K_2$. We claim that it has center $q$. Since homothecies send support real-hyperplanes (of $H\cap K_1$) to parallel support real-hyperplanes (of $H\cap K_2$), it has to send $H\cap \Delta$ to itself because both $K_1$ and $K_2$ are strictly convex and the homothecy is positive. So, it
must have center $q=K_1\cap\Delta=K_2\cap\Delta$.

Note that all these homothecies with center at $q$ (one for every hyperplane containing $L$) coincide at 
$L$ and  send $L\cap K_1$ to $L\cap K_2$.  Therefore, the ratio of homothecy coincides for every hyperplane $H$ through $L$. Then, because every point lies in one of the hyperplanes through $L$ (recall that $n\geq 3$), the  homothecy with center at $q$ and this ratio sends  $K_1$ to $K_2$ as we wished to prove.  \qed

As a corollary we have: 

\begin{cor}\label{corrg}
Let $K_1, K_2\subset \mathbb{C}^n$ be strictly convex bodies with $p_1\in int K_1$ and $p_2\in int K_2$, $n\geq 3$. Suppose that for every
hyperplane $H$ through the origin  
the section $(H+p_1)\cap K_1$ is a translated copy of $(H+p_2)\cap K_2$. Then, $K_1$ is a translated copy of $K_2.$ 
\end{cor}

\pf Translations are homothecies of ratio 1. So, by Theorem \ref{lemrg}, we conclude that there is a positive homothecy $h$ sending $K_1$ to $K_2$, with
ratio of homothecy $r>0$. To prove that  $h$ is a translation we must show that $r=1$.  Let $p_3=h(p_1)$. Then, for every hyperplane $H$ through the origin  
the section $(H+p_3)\cap K_2$ is  positively homothetic to $(H+p_2)\cap K_2$, with ratio of homothecy $r$. Let $H'$ be a hyperplane through the origin with the property that $(H'+p_3)=(H'+p_2)$. Hence, $(H'+p_3)\cap K_2=(H'+p_2)\cap K_2$, and being positively homothetic with ratio of homothecy $r>0$, we conclude that $r=1$. \qed

\begin{teo}\label{lemBB}
In dimension $n \geq 3$, a convex body $K$ all of whose hyperplane sections through an interior point $x_0$ are symmetric 
is itself symmetric. 
\end{teo}

\pf  To prove that $K$ is symmetric, it is enough, by Theorem \ref{teosc},  to prove that for every $\xi\in \mathbb{S}^1$, $\xi K$ is a translated copy of $K$. 

Fix  $\xi\in \mathbb{S}^1$, thus $\xi x_0$ is an interior point of $\xi K$.  It is enough to prove that for every hyperplane $H$ through the origin, $(H+x_0)\cap K$ is a translated copy of  $(H+\xi x_0) \cap \xi K$; because this implies, by Corollary \ref{corrg}, that $\xi K$ is a translated copy of $K$.

By hypothesis, 
$(H+x_0)\cap K$ is symmetric, hence by Theorem\ref{teosc}, $(H+x_0)\cap K$ is a translated copy of $\xi\big((H+x_0)\cap K\big)$, but 
$$\xi\big((H+x_0)\cap K\big)= (\xi H+\xi x_0)\cap \xi K= (H+\xi x_0)\cap \xi K.$$ \qed

Perhaps Theorem~\ref{lemBB} is false for $n=2$.

We now turn our attention to characterizing symmetry by means of projections.

The affine image of a symmetric set is a symmetric set, and the center is sent to the center. This is so, because if $\xi K=K$ for every $\xi\in \mathbb S^1$ and $f$ is a linear map, then $\xi f(K)=f(\xi K)=f(K)$. In particular, the orthogonal projection of a symmetric set onto any line is a disk; and this is also enough to have symmetry when the center is fixed:

\begin{lema}\label{lemSD}
A convex body $K$ is symmetric if and only if it has a translated copy all of whose orthogonal projections onto $1$-dimensional subspaces are disks centered at the origin.
\end{lema}

\pf The necessity is clear. For the sufficiency, assume $\pi(K)$ is a disk centered at the origin for every orthogonal projection $\pi$ onto a $1$-dimensional subspace.

We shall prove that for every $\xi\in \mathbb S^1$, $\xi K=K$.  

Suppose not. Then, there is $x\in K$ and $\xi\in \mathbb S^1$ such that 
$\xi x\notin K$. By convexity, there is a real $(2n-1)$-dimensional hyperplane $\Delta$ in $\mathbb R^{2n}=\mathbb C^n$ that separates $\xi x$ from $K$. Let $H$ by the unique hyperplane contained in $\Delta$. Let $\pi:\mathbb C^n \to L$ be the orthogonal projection onto the $1$-dimensional subspace $L$ orthogonal to $H$. Then the real line $\pi(\Delta)$ separates $\pi(K)$ from $\pi(\xi x)=\xi(\pi(x))$ on the line $L$. On the other hand $\pi(x)$ lies in the disk
$\pi(K)$ centered at the origin, therefore $\xi(\pi(x))$ lies in $\pi(K)$, which is a contradiction. \qed 

\begin{teo}\label{teoSim} 
For $k\geq 2$, a convex body all of whose orthogonal projections onto $k$-planes are symmetric is symmetric.
\end{teo}

\pf Consider a convex body $K\subset \mathbb C^n$, $n\geq 3$, all of whose orthogonal projections onto $(n-1)$-planes are symmetric. For every line $L$ through the origin, let $c_L$ be the center of the orthogonal projection of $K$ to $L^\perp$, the orthogonal hyperplane to $L$, and let $L^\prime$ be the parallel line to $L$ through $c_L$. We claim that
$$\bigcap_L L^\prime \text{ is a single point } \{x_0\}.$$

Let $L_1$ and $L_2$ be two different lines through the origin. The orthogonal projection of $K$ to $L_1^\perp\cap L_2^\perp$ is symmetric because it is the affine image of symmetric convex bodies (in $L_1^\perp$ and $L_2^\perp$); let $a$ be its center of symmetry and
let $\Gamma$ be the orthogonal plane to $L_1^\perp\cap L_2^\perp$ passing through $a$. The lines $L_1^\prime$ and $L_2^\prime$ are both in $\Gamma$ because the symmetry centers $c_{L_1}$ and $c_{L_2}$ are orthogonally projected to $a$. Therefore, since $L_1^\prime$ and $L_2^\prime$ are not parallel, they intersect in a point $\{x_0\}=L_1^\prime\cap L_2^\prime$.

Consider a third line $L_3$ through the origin, linearly independent to $L_1$ and $L_2$; it exists because $n\geq 3$. By the above argument, $L_3^\prime$ intersects both $L_1^\prime$ and $L_2^\prime$. But by the linear independence, $L_3^\prime$ intersects $\Gamma$ in at most one point, so we must have that $x_0\in L_3^\prime$.

Finally, for the general line $L$ through the origin, $L$ is linearly independent to at least one pair of the lines $L_1, L_2, L_3$. So that the preceding argument yields that $x_0\in L^\prime$. This proves that $\bigcap_L L^\prime =\{x_0\}$ as we claimed.

Suppose without loss of generality that $x_0$ is the origin. Then, for every $(n-1)$-plane $H$ through the origin, the orthogonal projection of $K$ 
onto $H$ is symmetric with center at the origin. This immediately implies that  orthogonal projections of $K$ onto $1$-dimensional subspaces are disks centered at the origin and therefore, by Lemma \ref{lemSD}, that $K$ is symmetric. 

The theorem now follows by induction, because a convex body all of whose orthogonal projections onto $k$-planes are symmetric has the property that all orthogonal projections onto $(k+1)$-planes are symmetric, and so on. 
\qed

The interesting question that remains open is then:  \emph{is a convex body all of whose orthogonal projections onto lines are disks, symmetric?}

\section{ Bombons are ellipsoids}

 Our main result is the following.
 
 \begin{teo}\label{teoBom}
A convex body $K\subset \mathbb{C}^n$ is an ellipsoid  if and only if for every line $L\subset \mathbb{C}^n$, $L\cap K$ is either empty, a point or a disk. 
\end{teo}

For expository reasons, it is convenient to define a \emph{bombon} as a convex body $K\subset \mathbb{C}^n$ such that for every line $L\subset \mathbb{C}^n$, $L\cap K$ is either empty or a disk; where points are regarded as disks of radius zero. Then, the theorem simply says that bombons and ellipsoids are one and the same.

Note that $1$-dimensional ellipsoids and $1$-dimensional  bombons are disks. Furthermore, complex ellipsoids are bombons because balls are bombons and the image of a bombon under an affine map is a bombon. The basic reason for this last fact, is that affine maps between lines send disks to disks. In \cite{ABM} (where bombons were considered as the complex case of what were called \emph{``pucks''}), it was proved that if a bombon has the property that the centers of all disks which are line sections parallel to a given line, lie in a hyperplane, then the bombon is an ellipsoid. Now, we shall prove that this hypothesis  is unnecessary.

The proof of Theorem~\ref{teoBom} is by induction. The general case will follow from a theorem proved in \cite{BM}, and the first case, $n=2$, takes most of the work here, using ideas of its own.

An \emph{abstract linear space} consists of a set $X$  together with a distinguished family of subsets, called \emph{abstract lines}, satisfying the following property:  given different $x,y\in X$, there is a unique abstract line $L$ containing $x$ and  $y$.  Typical examples of abstract linear spaces are euclidean $n$-space,  complex $n$-space and their associated projective spaces. 

Our interest lies in the abstract linear space $LS^3$, where $X=\mathbb S^3$ (the unit sphere of $\mathbb C^2$), and an abstract line is the intersection with $\mathbb S^3$ of a non-tangential complex line of  $\mathbb C^2$, which topologically is a flat circle in $\mathbb S^3$. 

A subset $A$ of an abstract linear space $X$ is \emph{linearly closed} if for any $x,y\in A$, the line through $x$ and $y$ is contained in $A$.  Since the intersection of linearly closed subsets is linearly closed, given any $Y\subset X$, there exists a unique minimal linearly closed subset containing $Y$, called its \emph{linear closure}.

\begin{lema}\label{lemclo}
Suppose $K\subset\mathbb{C}^2$ is a bombon with the property that the ellipsoid of minimal volume containing it, is the unit ball. Then, $K\cap \mathbb S^3$ is a linearly closed subset of the abstract linear space $LS^3$.
\end{lema}

\pf  Let $x,y\in K\cap \mathbb S^3$, $x\neq y$, and let $L$ be the line through $x$ and $y$. We have to prove that $L\cap \mathbb S^3\subset K$. 

Let $B$ be the unit ball in $\mathbb{C}^2$, whose boundary is $\mathbb S^3$. By hypothesis, $L\cap K$ and $L\cap B$ are two disks in $L$ which share two different points, $x$ and $y$ on their boundary. Then,  
$ K \subset B$ implies that $L\cap K = L\cap B$, and therefore that $L\cap \mathbb S^3\subset K$. \qed

\begin{prop}\label{prop:lin_cl}
The linearly closed proper subsets of the abstract linear space $LS^3$ are single points and lines.
\end{prop}

The proof of this proposition follows immediately from the next two lemmas. 

We need to make precise some standard definitions. 
A topological space $U$ can be topologically embedded in $\mathbb C$, if there is a continuous inyective map $f: U\to\mathbb C$. A topological space $A$ is locally embedded in $\mathbb C$ if for every point $x\in A$, there is a neighborhood $U$ of $x$ in $A$ which can be topologically embedded in $\mathbb C$.

\begin{lema}
Let $A$ be a linearly closed proper subset of the abstract linear space $LS^3$. Then $A$ is embedded in $\mathbb C$.
\end{lema}

\pf Since $A$ is proper, there exists a point $w \in\mathbb S^3$ such that $w\notin A$.  Let $L$ be the line through $w$ and the origin and let  $L^\perp$ be the orthogonal line to $L$ at the origin, which we identify with $\mathbb C$. 

Let $\pi:\mathbb S^3\setminus \{w\} \to L^\perp$  be the geometric projection from $w$, that is, given $x\in \mathbb S^3$ different  from $w$, $\pi(x)$ is the intersection with $L^\perp$ of the line through $x$ and $w$; well defined because this line is not parallel to  $L^\perp$.

We claim that $\pi |_A$ is an embedding. Indeed, if $x, y \in A$ are such that $\pi(x)=\pi(y)$ then their lines through $w$ coincide. If $x$ and $y$  were different points, it would imply that $w\in A$ because $A$ is linearly closed, which contradicts the choice of $w$.        \qed 

\begin{lema}
Let $A$ be a linearly closed, non-empty subset of the abstract linear space $LS^3$, which is not a point or a line. Then $A$ is not locally embedded in $\mathbb C$.
\end{lema}

\pf $A$ has at least two points, so it contains a line $\ell\subset A\subset\mathbb S^3$. But since it is not a line, there is yet another point $a\in A\setminus\ell$.

For every $x\in \ell$ let $\ell_x$ be the line through $a$ and $x$.  Note that $\ell_x\cap \ell=\{x\}$, and if $x\not=x'\in \ell$,
then $\ell_x\cap\ell_{x'}=\{a\}$, because two lines intersect in at most one point. Let 
$$\Omega=\bigcup_{x\in \ell}\ell_x \subset A\subset \mathbb S^3\,.$$  
Since every abstract line of $LS^3$ is actually a circle, then $\Omega\setminus\{a\}=\bigcup_{x\in \ell}(\ell_x\setminus\{a\})$ is an $\mathbb R$-bundle over $\mathbb S^1$. This implies that, either $\Omega$ is homeomorphic to the projective plane or $\Omega$ is homeomorphic to a closed cylinder modulo its boundary,
$\frac{\mathbb S^1\times [0,1]}{\mathbb S^1\times\{0,1\}}$, where the point $a\in \Omega$ is precisely $\frac{\mathbb S^1\times \{0,1\}}{\mathbb S^1\times\{0,1\}}$.

Since $\Omega\subset \mathbb S^3$, the space $\Omega$ cannot be homemorphic to the projective plane. So $\Omega$ is homeomorphic to $\frac{\mathbb S^1\times [0,1]}{\mathbb S^1\times\{0,1\}}$ and consequently, any neighborhood of $a$ in $\Omega$ cannot be topologically embedded in $\mathbb C$, because it contains the cone of two disjoint circles. Since $\Omega \subset A$, then $A$ is not locally embedded in $\mathbb C$.  \qed

\bigskip

\noindent \emph{Proof of Theorem~\ref{teoBom}.} As we have said, the proof is by induction on the dimension $n$. 

\smallskip
\noindent \emph{Case} $n=2$. Suppose $K\subset \mathbb{C}^2$ is a bombon, and let 
$E$ be the ellipsoid of minimal volume containing $K$. We may assume without loss of generality that $E$ is the unit ball of $\mathbb{C}^2$; and we must prove that $K = E$, which is equivalent to $K\cap \partial E=K\cap \mathbb S^3=\mathbb S^3$.

By Lemma \ref{lemclo},  $K\cap \mathbb S^3$ is a linearly closed subset of $LS^3$ and by Proposition \ref{prop:lin_cl}, we have to prove that $K\cap \mathbb S^3$ is neither a point nor a line. It clearly cannot be a point by the minimality condition. So, we are left to prove that this condition also rules out the line. We give an \emph{ad hoc} proof of this fact, without using John's theory.

Suppose there is a line $L$, such that $K\cap \mathbb S^3=L\cap \mathbb S^3$. We consider two cases. If $L$ passes through the origin, we may assume that
$L=\mathbb C\times\{0\}$. 

Consider the following family of ellipsoids for small $\epsilon\geq 0$, (see Fig. 1.a) 
$$E_\epsilon= \{(x,y)\in \mathbb{C}^2\mid \frac{|x|^2}{ (1+\epsilon)^2}+ \frac{|y|^2}{ (1-\epsilon)^2}\leq 1\}\,.$$
\begin{figure}[h!]
          \includegraphics[width=10cm]{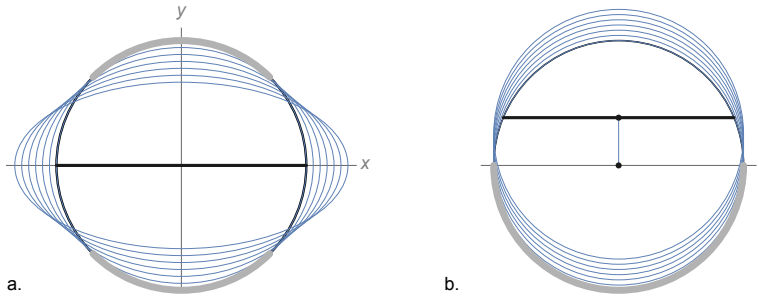}    
          \caption{\small {\bf a)} Section of the family of ellipsoids $E_\epsilon$, with a real plane having real lines in the coordinate axis.  {\bf b)} Section of the family of balls $E+\epsilon c$ with a real plane containing $c$ and the origin.}  
\end{figure}
 
For $\epsilon >0$, the volume of $E_\epsilon$ is less than the volume of the unit ball $E=E_0$; in fact, their volume ratio is $(1-\epsilon^2)^2$. Then, all the sets  $X_\epsilon=E\setminus E_\epsilon$ intersect $K$ because of the minimality condition of $E$. Furthermore, the intersection and limit of $X_\epsilon$, when $\epsilon \to 0$, is the solid torus $T=\{\,(x,y)\in\mathbb S^3\, |\, |y|\geq1/\sqrt{2}\,\}$. So that the compacity of $K$ gives us a point in $K\cap T\subset K\cap \mathbb S^3$, which is far from $L$.

If $L$ does not contain the origin, let $c$ be its closest point to the origin and consider the family of balls $E+\epsilon c$, $\epsilon\geq 0$ (see Fig.1.b). As before, the sets $X_\epsilon=E\setminus (E+\epsilon c)$ intersect $K$ for all $\epsilon>0$, because of the uniqueness (see Theorem 5 in \cite{ABM}) of the minimal volume ellipsoid containing $K$. Therefore, $K$ has a point in the limit of $X_\epsilon$, as $\epsilon \to 0$, which is the closed half-hemisphere of $\mathbb S^3$ opposite to $c$. Contradicting that $K\cap \mathbb S^3=L\cap \mathbb S^3$, and completing the proof for $n=2$.

\smallskip
\noindent\emph{General case.} Suppose the theorem is true for $n-1$, we shall prove it for $n\geq3$.  By induction, for every hyperplane $H$, the section $H\cap K$ is either empty, a single point or an ellipsoid. This implies that every section of $K$ is symmetric and consequently, by Theorem~\ref{lemBB}, $K$ is symmetric.  

We may assume without loss of generality that the center of $K$ is the origin.  Now note, by induction, that every hyperplane section of $K$ through the origin is an ellipsoid.
Then, Theorem~3.3 of \cite{BM}  implies that $K$ is an ellipsoid.  \qed

\section{Ellipsoids by means of sections and of projections}

Our main purpose now, is to prove that all sections through a point (orthogonal projections) of a convex body are ellipsoids only when then the body is an ellipsoid.  This result, under the additional hypothesis of symmetry, is Theorem 3.3 of \cite{BM}. We prove it here in general. 

\begin{teo}\label{teoksec}
A convex body in $\mathbb C^n$ all of whose $k$-sections through a point $p_0$  are ellipsoids is an ellipsoid, $2\leq k < n$. 
\end{teo}

\pf  It is easy to see that the general case $2\leq k < n$ of this theorem follows from the case $k=n-1$. 
 
Let $K\subset \mathbb{C}^n$ be a convex body with $n\geq 3$, and suppose  that all hyperplane sections of $K$ through $p_0$ are ellipsoids.
Let us first prove that $K$ is a real ellipsoid. Using Theorem 2.12.4 of \cite{MMO}, it is enough to prove that all real plane sections of $K$  are real ellipses.  Let $P$ be a real $2$-dimensional plane through $p_0$ and let $H$ be a hyperplane through $p_0$ such that $P\subset H\subset \mathbb{C}^n$. This is possible because $n\geq 3$. By hypothesis, $H\cap K$ is an ellipsoid.  Then, $P\cap K$  is a real ellipse.  This implies that $K$ is a real ellipsoid. 

Our next purpose is to prove that all  line sections of $K$ through $p_0$ are disks. Since every hyperplane section of $K$ through $p_0$ is an ellipsoid, we have that  all line sections of $K$ through $p_0$ are disks. Remember that any two parallel sections of a real ellipsoid are homothetic.  Hence, every line section, not only the line sections through $p_0$ are disks and consequently $K$ is a bombon.  By Theorem \ref{teoBom}, $K$ is an ellipsoid. 
\qed

\medskip
The analogous result for orthogonal projections is also true.  
 
\begin{teo}\label{teokproj}
A convex body in $\mathbb C^n$ all of whose orthogonal projections onto $k$-dimensional subspaces  
are ellipsoids is an ellipsoid, $2\leq k < n$. 
\end{teo}

\pf Suppose the orthogonal projection of the convex set $K\subset \mathbb C^n$ onto every $k$-dimensional subspace is an ellipsoid, $2\leq k < n$.
Since every ellipsoid is symmetric, by Theorem \ref{teoSim}, $K$ is symmetric. Furthermore, $K$ is a real ellipsoid. Let $P$ be a real $2$-dimensional plane and let $H$ be a $k$-subspace such that $P\subset H\subset \mathbb{C}^n$. This is possible because $k\geq 2$. By hypothesis, the orthogonal projection of $K$ onto $H$ is an ellipsoid $C$.  Then, the orthogonal projection of $K$ onto $P$ is a real ellipse.  By Theorem 2.12.5 of \cite{MMO}, $K$ is a real ellipsoid.  This concludes the proof because by Theorem 4.3 of \cite{BM}, every symmetric real ellipsoid is an ellipsoid. \qed

\medskip
The analogous over the real numbers of the following Theorem, holds only in dimensions greater than 2. It was originally proved in 1889 by H. Brunn, 
\cite{Bru}, under the hypothesis of regularity, and in its full generality  by G.R. Burton in 1977, \cite{B2}.
Our main result allows us to prove it for all dimensions in the complex case.

\begin{teo}\label{teocs}
A convex body all of whose hyperplane sections are symmetric is an ellipsoid. 
\end{teo}

\pf Let $K$ be a convex body all of whose hyperplane sections are symmetric and let  $L\subset \mathbb{C}^{n-1}$ be a line such that $L\cap K\not=\emptyset$. We shall prove that $L\cap K$ is  a disk and the theorem follows from Theorem~\ref{teoBom}. 
For that purpose, 
it will be enough to prove that there exists a hyperplane $H$ containing $L$,  such that the center of $H\cap K$ lies in $L$.  

Suppose, without loss of generality, that
$L=\{(0, \dots,0, z)\mid z\in \mathbb{C}\}$, and that the origin is in $L\cap K$.

For every $(n-2)$-plane  through the origin  $\Gamma\subset \mathbb{C}^{n-1}$, we have that 
$\Gamma\cap K\not=\emptyset$. Let $\Gamma^\prime$ be the hyperplane of $\mathbb{C}^n$ generated by $\Gamma$ and $L$, and let $x_\Gamma$ be the center of symmetry of $\Gamma^\prime\cap K$. If $x_\Gamma \in L$,
there is nothing to prove; so we may assume, without loss of generality, that for every $(n-2)$-plane  through the origin  $\Gamma$, we have that $x_\Gamma\notin L$.

Let $\pi: \mathbb{C}^{n}\to \mathbb{C}^{n-1}$ be the projection  onto the first $n-1$ coordinantes. The choice $\Gamma\to \pi(x_\Gamma)\in\Gamma$ is a continuous assignment of a non-zero vector in $\Gamma$, for each hyperplane $\Gamma$ of $\mathbb{C}^{n-1}$. This is a contradiction to the well known fact that the cannonical vector bundle of hyperplanes through the origin of $\mathbb{C}^{n-1}$ does not admit a non-zero section. See Steenrod's Book \cite{S}.

The continuity of $\Gamma\to \pi(x_\Gamma)$ follows from the following fact: suppose $H_i, i\geq 1$, is a sequence of hyperplanes in $\mathbb{C}^{n}$ with the property that 
$H_0=\lim 
(H_i)$ (where limits are taken as ${i\to\infty}$). If $H_i\cap K$ is symmetric with centre at $c_i$ for $i\geq 1$, then  $H_0\cap K$ is  symmetric with center at ${c_0=\lim c_i}$.
To see this, let $\Gamma_i=H_i-c_i$ and $K_i=(H_i\cap K)-c_i$, for $i\geq 0$.  Then, $\Gamma_0=\lim \Gamma_i$,  $K_0=\lim K_i$ (using the Hausdorff metric) and $K_i$ is symmetric with center at the origin for $i\geq 1$. This implies for every 
$\xi\in\mathbb{S}^1$ and $i\geq 1$,  that $\xi  K_i= K_i$. Therefore, 
$$\xi K_0=\xi(\lim K_i)=\lim \xi K_i =\lim  K_i= K_0$$ 
which implies that $H_0\cap K$ is symmetric with center at $c_0$.  \qed

 \medskip

The following, and final, result is again a characterization of ellipsoids by means of their line sections, but this time following the spirit of the Isometric Banach Conjecture (see  \cite{B},  \cite{BHJM}, and  \cite{BM}).  Our characterization deals with 1-dimensional sections that are not empty or a single point, which we call \emph{non-trivial}.

\begin{teo}\label{entenado}
 A convex body $K$ is an ellipsoid  if and only if any two  non-trivial line sections of $K$ are 
affinely equivalent. 
\end{teo}

\pf  In view of our main result, Theorem~\ref{teoBom}, it will suffice to prove that every non-trivial line section of $K$ is a disk. 

By hypothesis, there is a fixed convex set $C\subset \mathbb{C}$, such that every non-trivial line section of $K$  is affinely equivalent to $C$. Without loss of generality, we may assume that $\mathbb{S}^1$ is the circumcircle of $C$, so that
$G_C$, \emph{the compact group of  affine equivalences of $\mathbb{C}$ that fix $C$}, is a subgroup of $U(1)$. To see that every non-trivial line section of $K$ is a disk, we will prove that $C$ is the unit disk, which is clearly equivalent to  proving that $G_C=U(1)$.

Consider a point in the interior of $K$ and a 2-dimensional affine subspace through it. To fix ideas, we may assume that the origin is in the interior of $K$, and think of $\mathbb{C}^2\subset\mathbb{C}^n$.  The space of all lines through the origin in $\mathbb{C}^2$ is the Riemann sphere $\mathbb{S}^2$. Let 
$$\xi:\mathbb{\it E}\to \mathbb{S}^2$$
be the canonical vector bundle of $1$-dimensional subspaces of $\mathbb{C}^2$. 
In view of the continuity of the line sections of $K$ that we described in our previous proof, our hypothesis imply that 
the structure group of $\xi$ can be taken to be $G_C$. That is, if  $\xi$ is described by a characteristic map $\tau:\mathbb{S}^1\to U(1)$, it factors through $G_C$; or if it is described by  transition functions and an open cover, the transition functions can be chosen to go to $G_C$ (see, e.g., \cite{S}). 

Since the only proper compact subgroups of $U(1)$ (which can be identified with $\mathbb{S}^1$) are finite, but the canonical line bundle is non trivial, then $G_C=U(1)$. This completes the proof.
 \qed

\bigskip

\sn{\bf Acknowledgments.} Luis Montejano acknowledges  support  from CONACyT under 
project 166306 and  from PAPIIT-UNAM under project IN112614.

\end{document}